\documentclass{amsart}
\usepackage{amssymb, amscd}

\input xypic

\theoremstyle{plain}
\newtheorem{theorem}{Theorem}[section]
\newtheorem{corollary}[theorem]{Corollary}
\newtheorem{lemma}[theorem]{Lemma}
\newtheorem{proposition}[theorem]{Proposition}
\newtheorem{example}[theorem]{Example}
\theoremstyle{definition}
\newtheorem{definition}[theorem]{Definition}

\theoremstyle{remark}

\numberwithin{equation}{section}

\newcommand{\Res}{\operatorname{Res} }

\newcommand{\Proj}{\operatorname{Proj} }
\newcommand{\Sym}{\operatorname{Sym} }

\renewcommand{\S}{\Sigma}

\newcommand{\p}{\mathbf{P} }

\begin{document}

\title{Resultants and symmetric products}
\author{Helge Maakestad }
\address{Institut de Math{\`e}matiques, Universite Paris VII}
\email{maakestad@math.jussieu.fr}
\thanks{ }
\date{2003}
\begin{abstract} We use the symmetric product $\Sym^n(\p^1_k)$ of the projective line 
to describe the resultant scheme $R_{n,m}$ in $\p^n_k\times \p^m_k$ as a quotient $X/G$ 
where $X=(\p^1_k)^{n+m}$ and $G\subseteq Aut_k(X)$ is a finite subgroup. As a special case
we give a description of the discriminant scheme in terms of the symmetric product.
\end{abstract}
\maketitle

\tableofcontents
\section*{Introduction} In this paper we will consider the resultant scheme $R_{n,m}$ whose closed points parametrize
pairs of homogeneous polynomials $(f,g)$ of degree $n$ and $m$ in two variables with common roots. 
Given general degree $n$ and $m$ homogeneous polynomials
$f(x_0,x_1)$ and $g(x_0,x_1)$ in two independent 
variables $x_0,x_1$ over a fixed algebraically closed basefield $k$ of characteristic
zero. There exists a
polynomial $\Res(f,g)$ in the coefficients of $f$ and $g$ which is zero if and only if $f$ and $g$ have a common root.
The polynomial $\Res(f,g)$ is an irreducible  bi-homogeneous polynomial defined up to a constant.
The \emph{resultant scheme of degree $n$ and $m$} denoted $R_{n,m}$ is by definition the projective 
subscheme of $\p^n_k\times \p^m_k$ defined by $\Res(f,g)$. This scheme has been studied by several authors
(see \cite{gkz}, \cite{jou}). 
The novelty of this paper is the geometric construction of $R_{n,m}$ using the symmetric product.

It is a standard fact that projective $n$-space might be described as follows: There
exists an action of the symmetric group $\S_n$ on $n$ letters on $n$ copies of the projective line $\p^1_k$: 
\[ \S_n\times (\p^1_k)^n \rightarrow (\p^1_k)^n ,\]
with the property that the quotient - the symmetric product - 
\[ \Sym^n(\p^1_k):=(\p^1_k)^n/\S_n \] 
is isomorphic to projective $n$-space. 

The \emph{resultant scheme} $R_{n,m}$ of two polynomials $f$ and $g$ is the  subscheme of $\p^n_k\times \p^m_k$ 
defined as the zero-set of the 
resultant polynomial $\Res(f,g)$. 
We prove in Theorem \ref{resultant} that the resultant scheme $R_{n,m}$ may be described as a quotient 
$X_{n,m}/\S_n\times \S_m$, where $X_{n,m}$ is a union of products of projective lines.

As a special case we consider the discriminant scheme $D_n$ which 
is a hypersurface in $\p^n_k$ and consider the scheme-theoretic
inverse image $X_n:=\pi^{-1}(D_n)$, where $\pi:(\p^1_k)^n\rightarrow \Sym^n(\p^1_k)$ is the canonical quotient-map.
We prove in Theorem \ref{Disc} that the discriminant scheme $D_n$ may be described as $X_n/\S_n$ where the $\S_n$-action on $X_n$
is the one induced from the $\S_n$-action on $(\p^1_k)^n$. Hence we get a quotient map
\[ X_n\rightarrow D_n ,\]
where $X_n$ is a union of products of projective lines  with non-reduced structure.

 It would be interesting to give a geometric construction of Jouanolou's multivariate resultant as found in 
\cite{jou} generalizing the construction in this paper, and I hope to return to this problem sometime in the future. 

\section{Elementary homogeneous symmetric polynomials}
In this section we define the \emph{elementary homogeneous symmetric polynomials}
and relate them to  the theory of the resultant of two
homogeneous polynomials in two variables, and the discriminant of a
homogeneous polynomial in two variables following \cite{LANG}. 

In general, let $A$ be a
commutative ring with unity, and let $Z,v_0,\cdots ,v_n$ 
and $w_0,\cdots ,w_m$ 
be algebraically independent variables over $A$. We consider two polynomials
\[ f_v(Z)=v_0Z^n+v_1Z^{n-1}+\cdots +v_n \]
\[ g_w(Z)=w_0Z^m+v_1Z^{m-1}+\cdots +w_m \]
and make the following definition:

\begin{definition} \label{DRes}The resultant of $f_v$ and $g_w$ denoted by
$\Res (f_v,g_w)$ is given by the following determinant
$$
\begin{vmatrix} 
v_0 & v_1 & v_2 & \cdots & v_{n-1} & v_n & 0 & \cdots  & 0  & 0 \\
0   & v_0 & v_1 & \cdots & v_{n-2} & v_{n-1} & v_n & \cdots & 0 & 0  \\
\vdots & \ddots&  & & & &\ddots  & \vdots& v_n  & 0 \\
0   & 0   & 0   & \cdots & 0 & v_0     & v_1     & v_2 & \cdots & v_n \\
w_1 & w_2 & w_2 & \cdots & w_{m-1} & w_m & 0 &  \cdots & 0 & 0\\
0 & w_1 & w_2 & \cdots & w_{m-2} & w_{m-1} & w_m & \cdots & 0 & 0
\\
\vdots &\ddots  & & & &\ddots &  \vdots & w_{m-1} & w_m & 0 \\ 
0 & 0   & 0    & \cdots & w_1   & w_2 & w_3  & \cdots & 0 & w_m \\
\end{vmatrix}
$$
\end{definition}

The resultant $\Res (f_v,g_w)$ is bihomogeneous of degree $(m,n)$ in the
variables $v$ and $w$. This follows directly from the definition. 

Let $x_{01},x_{11},\cdots ,x_{0n},x_{1n},y_{01},y_{11},\cdots ,
y_{0m},y_{1m},X,Y$ be algebraically independent over $A$.  Consider
the polynomials
\[ F_n(x_{ij},X,Y)=(x_{01}X-x_{11}Y)\cdots (x_{0n}X-x_{1n}Y)= \]
\[ p_0(x_{ij})X^n-p_1(x_{ij})X^{n-1}Y+\cdots +(-1)^np_n(x_{ij})Y^n \]
and
\[G_m(y_{ij},X,Y)=(y_{01}X-y_{11}Y)\cdots (y_{0m}X-y_{1m}Y) =\]
\[ p_0(y_{ij})X^m-p_1(y_{ij})X^{m-1}Y+\cdots +(-1)^mp_m(y_{ij})Y^m .\]
We make the following definition:
\begin{definition} \label{EHSP}The polynomials $p_0(x_{ij}),\cdots ,p_n(x_{ij})$
  are the \emph{elementary homogeneous symmetric polynomials} with respect
to the variables $x_{ij}$.
\end{definition}

\begin{example} We compute the elementary homogeneous symmetric
polynomials in an example.
\end{example}
Consider the case $n=2$: We get
\[ F_2(x_{ij},X,Y)=(x_{01}X-x_{11}Y)(x_{02}X-x_{12}Y)= \]
\[ x_{01}x_{02}X^2-(x_{01}x_{12}+x_{11}x_{02})XY+x_{11}x_{12}Y^2 \]
hence
\[ p_0(x_{ij})=x_{01}x_{02},\]
\[p_1(x_{ij})= x_{01}x_{12}+x_{11}x_{02},\]
and
\[ p_2(x_{ij})= x_{11}x_{12}.\]
\begin{proposition}\label{IND} The elementary homogeneous symmetric polynomials 
are algebraically independent over $A$.
\end{proposition}
\begin{proof} Consider the polynomial $F_n(x_{ij},X,Y)$, letting
$Y=x_{0i}=1$ for $i=1,\cdots ,n$, we get the polynomial
\[ (X-x_{11})\cdots (X-x_{1n})=X^n-s_1(x_{1i})X^{n-1}+\cdots
+(-1)^ns_n(x_{1i}). \]
Assume we have a non-trivial algebraic relation involving the
elementary homogeneous symmetric polynomials, letting $Y=x_{0i}=1$ for 
$i=1,\cdots ,n$, we get a nontrivial algebraic relation involving
the elementary symmetric polynomials, which we know are algebraically 
independent over $A$. This is a contradiction, hence we have proved the assertion.
\end{proof}

Put $Z=X/Y$ and $t_i=x_{1i}/x_{0i}, u_i=y_{1i}/y_{0i}$. We may
write
\[ F_n(x_{ij},X,Y)=x_{01}\cdots x_{0n}Y^n(Z-t_1)\cdots (Z-t_n)=\]
\[ Y^n(p_0(x_{ij})Z^n-p_0(x_{ij})s_1(t)Z^{n-1}+\cdots
+(-1)^np_0(x_{ij})s_n(t)) \]
and 
\[ G_m(y_{ij},X,Y)=y_{01}\cdots y_{0m}Y^m(Z-u_1)\cdots (Z-u_m) =\]
\[ Y^m(p_0(y_{ij})Z^m-p_0(y_{ij})s_1(u)Z^{m-1}+\cdots
+(-1)^mp_0(y_{ij})s_m(u)). \]
\textbf{Notation}: $v_0=p_0(x_{ij}), v_k=(-1)^kp_0(x_{ij})s_k(t)$ for 
$k=1,\cdots ,n$, and $w_0=p_0(y_{ij}),w_l=(-1)^lp_0(y_{ij})s_l(u)$ for
$l=1,\cdots ,m$. From proposition \ref{IND} it follows that the new variables $v_i,w_j$ are
algebraically independent over $A$. It also follows that the variables
$v_0,t_1,\cdots ,t_n,w_0,u_1,\cdots ,u_m$ are algebraically
independent over $A$, hence we may consider the following
factorization:

\[ F_n(x_{ij},X,Y)=Y^n(v_0Z^n+\cdots +v_n)=Y^nf_v(Z) \]
\[ G_m(y_{ij},X,Y)=Y^m(w_0Z^m+\cdots +w_m)=Y^mg_w(Z) .\]

We get a theorem: 

\begin{theorem} \label{Resth} With the previous definitions the following holds: 
\[ \Res (f_v,g_w)=\prod_{i=1}^n\prod_{j=1}^m(x_{1i}y_{0j}-x_{0i}y_{1j})=
v_0^mw_0^n\prod_{i=1}^{n}\prod_{j=1}^m(t_i-u_j). \]
\end{theorem}
\begin{proof} See \cite{LANG}, section VI.8.
\end{proof}

The formula in Theorem \ref{Resth} is related to the \emph{Poisson formula} as found in \cite{jou}. 

Note that if $f(Z),g(Z)$ are polynomials with coefficients in an algebraically closed field
it follows that $\Res (f,g)=0$ if and only if $f$ and $g$ have at least one common root.
We also see that $\Res (f,f')=0$ if and only if $f$ has a root of multiplicity at least 2.

\begin{definition} Given a polynomial 
$f_v(Z)=v_0(Z-t_1)(Z-t_2)\cdots (Z-t_n)=$ $v_0Z^n+\cdots v_n$
we define its discriminant $D(f_v)$ as follows:
\[ D(f_v)=(-1)^{\frac{n(n-1)}{2} }v_0^{2n-2}\prod_{i\neq j}(t_i-t_j)
\]
\end{definition}

\begin{lemma}\label{FactDisc}  $D(f_v)=v_0^{2n-2}\prod_{i<j}(t_i-t_j)^2=\prod_{i<j}
(x_{1i}x_{0j}-x_{0i}x_{1j})^2$
\end{lemma}
\begin{proof} Easy calculation. 
\end{proof}

We have a relation with the resultant:

\begin{proposition} \label{ResDisc}$\Res (f_v,f'_v)= v_0^{2n-1}\prod_{i\neq
    j}(t_i-t_j)= $
\[ (-1)^{\frac{n(n-1)}{2} }v_0D(f_v)=(-1)^{\frac{n(n-1)}{2} }
    v_0^{2n-1}\prod_{i<j}(t_i-t_j)^2 \]
hence formally we may write
\[D(f_v)=(-1)^{\frac{n(n-1)}{2}}\frac{1}{v_0}\Res (f_v,f_v').\]

\end{proposition}
\begin{proof} Easy calculation.
\end{proof}

Fix an algebraically closed field $k$ of characteristic zero. The set of homogeneous polynomials
$F_p(X,Y)=p_0X^n+p_1X^{n-1}Y+\cdots + p_nY^n$ in $k[X,Y]$ of degree $n$
is in one-to-one correspondence with the closed points of $\p^n_k$, with homogeneous coordinates $p_0,\cdots
,p_n$. The \emph{discriminant scheme of degree n binary forms} $D_n$ is by definition the 
closed sub-scheme of $\p^n_k$ defined by the discriminant polynomial 
$D(f_p)=\frac{1}{p_0}\Res (f_p,f_p')$, where
\[ f_p(Z)=p_0Z^n+p_1Z^{n-1}+\cdots +p_n.\]

\begin{lemma} (Sylvester formula) \label{syl} The polynomial $\Delta_n:=\frac{1}{p_0}\Res (f_p,f_p')$ is given by 
the following determinant
$$
\frac{1}{p_0}\begin{vmatrix} 
p_0 & p_1 & p_2 & \cdots & p_{n-1} & p_n & 0 & \cdots  & 0  \\
0   & p_0 & p_1 & \cdots & p_{n-2} & p_{n-1} & p_n & \cdots & 0  \\
\vdots &  &  & \ddots &\vdots  & \vdots  &   & \ddots & \vdots \\
0   & 0   & 0   & \cdots & p_0     & p_1     & p_2 & \cdots & p_n \\
p_1 & 2p_2 & 3p_2 & \cdots & (n-1)p_{n-1} & np_n & 0 &  \cdots & 0 \\
0 & p_1 & 2p_2 & \cdots & (n-2)p_{n-2} & (n-1)p_{n-1} & np_n & \cdots & 0
\\
\vdots &   & &\ddots  &\vdots  & \vdots &  & \ddots  & \vdots \\ 
0 & 0   & 0    & \cdots & p_1   & 2p_2 & 3p_3  & \cdots & np_n \\
\end{vmatrix}
$$
\end{lemma}
\begin{proof} This follows from proposition \ref{ResDisc} and
  definition \ref{DRes}.
\end{proof}

\begin{example}We compute the polynomial $\frac{1}{p_0}\Res
  (f_p,f'_p)$ in two examples.
\end{example}
 Using the Sylvester-formula \ref{syl} one computes  the
  discriminant of a quadratic form:

\[ \Delta_2=4p_0p_2-p_1^2 \]
and a cubic form:
\[ \Delta_3= 27 p_0^2p_3^2 + 4p_0p_2^3 + 4p_1^3p_3 - p_1^2p_2^2 - 18p_0p_1p_2p_3. \]

\section{Resultants and symmetric products}

In this section we consider the action of $\S_n$ on 
$\p^1_k\times \cdots \times \p^1_k$. We use the symmetric product $\Sym^n(\p^1_k)$ of the projective line
to describe the resultant-scheme $R_{n,m}$ and discriminant-scheme $D_n$ as quotients 
$X/G$  where $G\subseteq Aut_k(X)$ is a finite subgroup. Let in this section $k$ be an algebraically 
closed field of characteristic zero.

 Recall first some standard facts from invariant-theory for finite groups. Assume $Y$ is a projective sub-scheme of
$\p^n_k$ and $G$ a finite
group acting linearly on the homogeneous coordinate-ring $S$ of $Y$. We assume the action 
$\tilde{\sigma}:G\times S\rightarrow S$ 
is graded, ie $G(S_d)\subseteq S_d$, hence the invariant-ring $S^G$ is a graded subring of $S$.
The action of $G$ on $S$ induces an action

 \[  \sigma:G\times Y\rightarrow Y  \]

 and a quotient $Y/G$ exists. It is given by the inclusion of graded rings
$S^G \rightarrow S$.
The quotient-map for the action $\sigma$ is given by the natural map

\begin{gather}
 \pi:Y=\Proj(S)\rightarrow \Proj(S^G)=Y/G.\label{invariantring} 
\end{gather}
The map $\pi$ has the following properties: it is surjective, affine, maps $G$-invariant closed subsets
to closed subsets and separates $G$-invariant closed subsets. The orbits of the action $\sigma$ are closed
(see  \cite{SPR}).
We state a well known result.

\begin{proposition} \label{closed}
Let $Z$ be a closed subscheme of $Y/G$, and let $X=\pi^{-1}(Z)$ be the scheme-theoretic inverse
image. There exists an induced action of $G$ on $X$, and $X/G\cong Z$.
\end{proposition}
\begin{proof} This follows from Theorem 7.1.4, \cite{B}.
\end{proof}

We state another result which will be needed in this section:
Assume $W=\Proj(T)\subseteq \p^m_k$ is a projective scheme, and $\tilde{\tau}:H\times T\rightarrow T$
is a graded action of a finite group $H$ on the homogeneous coordinate-ring of $W$. 
The action $\tilde{\tau}$ induces an action $\tau:H\times W\rightarrow W$, and we get a product-action
\[ \rho:G\times H\times Y\times W\rightarrow Y\times W .\]

\begin{proposition}\label{product}  There exists an isomorphism
\[ Y\times W/G\times H \cong Y/G\times W/H .\]
\end{proposition}
\begin{proof} The proposition follows from Remark 7.1.7 \cite{B} as
  follows: Let $e$ denote the group consisting of the identity. We have an inclusion of groups 
\[ G\times e \subseteq G\times H \]
hence we get isomorphisms
\[ Y\times W/G\times H \cong Y\times W/ G\times e/e\times H \cong \]
\[ Y/G\times W/ e\times H \cong Y/G\times W/H, \]
and the Proposition follows. 

\end{proof}

 Let 
$X,Y,x_{01},x_{11},\cdots ,x_{0n},x_{1n}$ be algebraically independent variables
over $k$. Consider the polynomial
\[ F_n(x_{ij},X,Y)=(x_{01}X-x_{11}Y)\cdots (x_{0n}X-x_{1n}Y) \]
in $k[x_{ij},X,Y]$. Multiply the parentheses in the expression
\[ (x_{01}X-x_{11}Y)\cdots (x_{0n}X-x_{1n}Y) \]
to get a polynomial
\[ p_0(x_{ij})X^n-p_1(x_{ij})X^{n-1}Y+\cdots +(-1)^np_n(x_{ij})Y^n .\]
Recall from definition \ref{EHSP} that the polynomials
$p_k(x_{ij})$ are the elementary homogeneous symmetric polynomials with
respect to the variables $x_{ij}$.

Consider the map
\[ \pi: \p^1_k\times \cdots \times \p^1_k\rightarrow \p^n_k \]
defined on closed points by
\[ \pi((x_{01}:x_{11}),\cdots ,(x_{0n}:x_{1n}))=(p_0(x_{ij});\cdots 
; p_n(x_{ij}) ).\]
There exist an obvious action $\sigma$ of the symmetric group
$\S_n$ on $\p^1_k\times \cdots \times \p^1_k$, and the map $\pi$ 
is invariant with respect to $\sigma$ because of the following:
Consider a closed point $((x_{01}:x_{11}),\cdots ,(x_{0n}:x_{1n}))$ in $(\p^1_k)^n$. We may
consider the homogeneous polynomial
\[ (x_{01}X+x_{11}Y)\cdots (x_{0n}X+x_{1n}Y)=F_n(x_{ij},X,Y) .\]
Given a permutation $\sigma$ in $\S_n$, it acts naturally on $(\p^1_k)^n$ 
and permuting the linear terms in the polynomial $F_n(x_{ij},X,Y)$ corresponds 
to the action of $\S_n$ on $(\p^1_k)^n$. Permuting the linear terms in the polynomial $F_n(x_{ij},X,Y)$
does not change the polynomial $F_n(x_{ij},X,Y)$, and it follows that $\pi$ is $\S_n$-invariant.
By definition
\[ \p^1_k\times \cdots \times \p^1_k=\Proj k[x_{p(1)1}\cdots x_{p(n)n} ]
\]
where $p(i)=0,1$ for $i=1,\cdots ,n$. The ring $k[x_{p(1)1}\cdots
x_{p(n)n} ]$ is by definition the n-fold Cartesian product of graded rings
\[ k[x_{01},x_{11}]\times_k \cdots \times_k k[x_{0n},x_{1n}] .\]
(See \cite{HAR}, ex.II.5.11for the definition of Cartesian product).
The map $\pi:\p^1_k\times \cdots \times \p^1_k\rightarrow \p^n_k$
is defined by the map of graded rings
\begin{equation} \label{segre}
 \phi:k[y_0,\cdots ,y_n]\rightarrow k[x_{p(1)1}\cdots x_{p(n)n} ] 
\end{equation} 
where
\[ \phi(y_0) = p_0, \cdots ,\phi(y_n)=p_n.  \]
The action of $\S_n$ on $\p^1_k\times \cdots \times \p^1_k$
lifts to a linear graded action of $\S_n$ on the homogeneous coordinate-ring
$k[x_{p(1)1}\cdots x_{p(n)n} ]= k[x_{01},x_{11}]\times_k \cdots \times_k k[x_{0n},x_{1n}]$.
The action of $\S_n$ on $k[x_{01},x_{11}]\times_k \cdots \times_k k[x_{0n},x_{1n}]$ is the obvious one.

\begin{theorem} \label{INV}The invariant-ring $k[x_{p(1)1}\cdots x_{p(n)n}
  ]^{\S_n}$
is generated by the elementary homogeneous symmetric polynomials.
\end{theorem}
\begin{proof} This follows directly from \cite{LANG}, Theorem IV.6.1.

\end{proof}

Note that Theorem \ref{INV} is valid over any commutative ring $A$ with unit.

Consider the resultant polynomial $\Res(f_v,g_w)$ and the resultant scheme
$R_{n,m}$ defined as the zero-set of the resultant polynomial. 
We saw in section 1 that the resultant $\Res (f_v,g_w)$ of two polynomials
\[ f_v(X)=v_0X^n+\cdots +v_n \]
and
\[ g_w(X)=w_0X^m+\cdots +w_m,\]
is a bihomogeneous polynomial of degree $(n,m)$, hence it defines
a closed sub-scheme $R_{n,m}$ of $\p^n_k\times \p^m_k$. Here we let
$\p^n_k\times \p^m_k$ have homogeneous coordinates $v_0,\cdots ,v_n,w_0,\cdots ,w_m$.
We have a product-mapping
\[\pi_{n,m}:=\pi_n\times \pi_m: (\p^1_k)^n\times (\p^1_k)^m \rightarrow \p^n_k\times \p^m_k .\]
There is an obvious group-action 
\[ \sigma: (\S_n\times \S_m)\times  (\p^1_k)^n\times (\p^1_k)^m \rightarrow
(\p^1_k)^n\times (\p^1_k)^m, \]
and the pair $(\p^n_k\times \p^m_k,\pi_n\times \pi_m)$ is a 
quotient for the action $\sigma$ by proposition \ref{product}. Define 
$X_{n,m}$ as the scheme-theoretic inverse image $\pi_{n,m}^{-1}(R_{n,m})$. 
It follows that $X_{n,m}$ is a closed subscheme of
$(\p^1_k)^n\times (\p^1_k)^m $. The map  $\pi_{n,m}$ is $\S_n\times \S_m$-invariant hence there exists an action 
of $\S_n\times \S_m$ on $X_{n,m}$. We get an induced map 
\[ \tilde{\pi}_{n,m}:X_{n,m}\rightarrow R_{n,m}.\]

\begin{lemma} $X_{n,m}=V(\prod_{i,j}(x_{0i}y_{1j}-x_{1i}y_{0j}))$.
\end{lemma}
\begin{proof} The lemma follows from theorem \ref{Resth}.
\end{proof}

\begin{theorem} \label{resultant} The induced map $\pi_{n,m}: X_{n,m}\rightarrow R_{n,m}$
induces an isomorphism 
\[ X_{n,m}/(\S_n\times \S_m) \cong R_{n,m}.\]
\end{theorem}
\begin{proof} By proposition \ref{product} the map
\[\pi_{n,m}: (\p^1_k)^n\times (\p^1_k)^m \rightarrow \p^n_k\times
\p^m_k \]
is a quotient-map for the $\S_n\times \S_m$-action, hence by proposition \ref{closed} the induced map
\[ X_{n,m}\rightarrow R_{n,m} \]
is a quotient-map  for the $\S_n\times \S_m$-action on $X_{n,m}$.
\end{proof}
 
As a special case we consider the discriminant scheme of degree n binary forms $D_n$. 

\begin{corollary}\label{quotient} Consider the map 
\[ \pi: \p^1_k\times \cdots \times \p^1_k\rightarrow \p^n_k \]
given by
\[ \pi((x_{01}:x_{11}),\cdots ,(x_{0n}:x_{1n}))=(p_0(x_{ij});\cdots 
; p_n(x_{ij}) ),\]
where $p_0,\cdots ,p_n$ are the elementary homogeneous symmetric polynomials.
The map  $\pi$ is a quotient-map for the natural $\S_n$-action on 
$\p^1_k\times \cdots \times \p^1_k$ and it follows that 
$\Sym^n(\p^1_k)$ is isomorphic to $\p^n_k$. 
\end{corollary}
\begin{proof} By  remark \ref{invariantring} the quotient $(\p^1_k)^n/\S_n$ is given by
the map 
\[ \p^1_k\times \cdots \times \p^1_k=\Proj k[x_{p(1)1}, \cdots ,x_{p(n)n}]\rightarrow \Proj k[x_{p(1)1}
\cdots x_{p(n)n}]^{\S_n}. \]
By theorem \ref{INV} we have $k[x_{p(1)1}
\cdots x_{p(n)n}]^{\S_n}$ equals $k[p_0,\cdots ,p_n]$ where $p_0,\cdots ,p_n$ are the elementary homogeneous
symmetric polynomials. By proposition \ref{IND}
the elementary homogeneous symmetric polynomials are algebraically
independent over $k$, hence $\Proj k[p_0,\cdots ,p_n]=\p^n_k$,
and the result follows.
\end{proof}

Note that corollary \ref{quotient} gives an explicit quotient map  
\[ \pi:(\p^1_k)^n\rightarrow \Sym^n(\p^1_k)\cong \p^n_k .\]
It is defined in terms of the elementary homogeneous symmetric polynomials.
The map $\pi$ is studied in \cite{katz} in a similar
situation. It is referred to as the \emph{Viete map}.

The discriminant $D_n$ is a closed subscheme of $\p^n_k$ defined by the 
polynomial $\Delta_n$, and it is well known that $\Delta_n$ is an irreducible polynomial. 
It follows that $D_n$ is an irreducible hypersurface in
$\p^n_k$. Let $X_n$ be the scheme-theoretic inverse image $\pi^{-1}(D_n)$.
There exists an induced map $\tilde{\pi}:X_n\rightarrow D_n$. 
Since $\pi$ is $\S_n$-invariant there exist an induced action
of $\S_n$ on $X_n$ and the map $\tilde{\pi}$ is $\S_n$-invariant.  
We get a description of the scheme $X_n$:

\begin{lemma} \label{fdisc} $X_n=V(\prod_{i<j}(x_{1i}x_{0j}-x_{0i}x_{1j})^2)$
\end{lemma}
\begin{proof} This follows from lemma \ref{FactDisc}.
\end{proof}
 
Note that the irreducible components of $X_n$ are products of projective spaces with non-reduced structure.

\begin{theorem} \label{Disc} The induced map $\tilde{\pi}:X_n\rightarrow D_n$ is a
quotient-map for the natural $\S_n$-action on $X_n$, hence
we get an isomorphim 
\[ X_n/\S_n\cong D_n .\]
\end{theorem}
\begin{proof} Since $\pi:\p^1_k\times \cdots \times \p^1_k\rightarrow
  \p^n_k$ is a quotient-map, the theorem follows from proposition \ref{closed}.
\end{proof}

\textbf{Acknowledgements.}  The paper was revised during a stay at Universite Paris VII
financed by a fellowship from the RTN network HPRN-CT-2002-00287, 
Algebraic K-theory, Linear algebraic groups and related structures and
I would like to thank Max Karoubi for an inviation to Paris
VII. Thanks also to the referee and Ragni Piene for useful comments.


\begin{thebibliography}{4}
\bibitem{AM} M. F. Atiyah, Introduction to commutative algebra,
  \emph{Addison-Wesley Publishing Co.} (1969)

\bibitem{B} A. Bialynicki-Birula, Quotients by actions of groups,
  Encyclopediq Math. Sci. 131, \emph{Springer Verlag} (2002), 1-82


\bibitem{gkz} I. Gelfand, M. Kapranov, A. Zelevinsky, Discriminants, resultants and multidimensional 
determinants, \emph{Birkh{\"a}user Boston} (1994) 

\bibitem{HAR} R. Hartshorne, Algebraic geometry, GTM no. 52, \emph{Springer-Verlag} (1977)

\bibitem{jou} J.P. Jouanolou, Le formalisme du resultant,
  \emph{Adv. Math. } (1991), 90, 117-263

\bibitem{katz} G. Katz, "How tangents solve algebraic equations, or a
  remarkable geometry of discriminant varieties", \emph{Expo. Math.}
  (2003) 3, 219-261


\bibitem{LANG} S. Lang, Algebra, \emph{Addison-Wesley} (1993)

\bibitem{MFK} D. Mumford, Geometric invariant theory, \emph{Springer-Verlag} (1991)

\bibitem{SPR} T. A. Springer, Invariant theory, Lecture Notes in Mathematics vol. 585,
\emph{Springer-Verlag} (1977)
\end{thebibliography}
\end{document}